\documentclass[12pt]{article}
\usepackage[pdftex,letterpaper=true,colorlinks=true,pdfpagemode=none,urlcolor=blue,linkcolor=blue,citecolor=blue,pdfstartview=FitH]{hyperref}
\usepackage{amsmath,amsfonts}
\usepackage{graphicx}
\usepackage{color}
\usepackage[labelformat=simple]{subcaption}
\usepackage{cases}
\usepackage{authblk}

\usepackage{tikz}
\title{Mean conservation for density estimation via diffusion using the finite element method}
\author{Keith Y. Patarroyo\thanks{kypatarroyot@unal.edu.co, Undergraduate physics student at Universidad Nacional de Colombia, sede Bogot\'a. }}

\begin{document}

\maketitle
\begin{abstract}

We propose  boundary conditions for the diffusion equation that maintain the initial mean and the total mass of a discrete data sample in the density estimation process. A complete study of this framework with numerical experiments using the finite element method is presented for the one dimensional diffusion equation, some possible applications of this results  are presented as well. We also comment on a similar methodology for the two-dimensional diffusion equation for future applications in two-dimensional domains.
\end{abstract}

\section{Introduction}

Estimating a density function using a set of initial data points in order to find probability information is a very significant tool in statistics\cite{bib1}. The method of Kernel Density Estimation (KDE)\cite{bib2} is now standard in many analysis and applications. Furthermore, this idea  has been applied in multiple fields (Archaeology \cite{bib3}, Economy \cite{bib4}, etc). The author of this article is particularly interested in constructing Perception of Security (PoS) hotspots using (KDE) methods to analyze real data registered by security experts in Bogot\'a \cite{bib0}. 

Nowadays a wide variety of methods are available to find density functions \linebreak(KDE) \cite{bib1},\cite{bib6}. The  method of KDE via difussion is of particular interest for this document; a recent article  \cite{bib5} develops a systematic method for (KDE) using the diffusion equation, also they propose a more general equation to solve some biases for data estimation. However in their analysis, it is only considered the normalization (conservation of mass) of the density function via Neumann boundary conditions, the mean of the sample data is not considered, thus inducing a change of an important initial parameter  from the discrete data sample.
 
In this article, we propose a new set of boundary conditions for the diffusion equation that maintain the initial mean and mass of the the discrete data sample in the density estimation process. A complete study of this framework is performed  using the finite element method (FEM) to solve the one-dimensional diffusion equation for different boundary conditions. We show the induced error on the final density when the mean is not conserved. We also show how this one-dimensional model can be used to simulate a (PoS) in a busy avenue of a city. Lastly the new boundary conditions are presented for the two-dimensional diffusion equation for future applications in two dimensional domains. 

\section{Diffusion equation with different boundary conditions}

As it was first noted in \cite{bib6} and expanded in \cite{bib5}, solving the diffusion equation with a discrete data sample $\{b_n\}_{n=1}^{N}$ as initial condition (\ref{eq2})  give an estimate of a continuous probability density function. Then by solving the diffusion equation \cite{bib10},

\begin{numcases}{}
  \frac{\partial u(x,t)}{\partial t}- \frac{\partial^2 u(x,t)}{\partial x^2}=0 \qquad a<x<b , t>0, \label{eq1}
   \\
     u(x,0)=\frac{1}{N}\sum_{i=1}^{N}\delta(x-b_i), \ \ x,b_i  \in [a,b] , \label{eq2}
\end{numcases}

with appropriate boundary conditions and then finding the best $t$ (bandwidth) for the initial data sample one obtains a continuous estimation of the experimental density. In this article we do not consider algorithms for bandwidth selection, we consider only the conservation of the mean. For more information on the bandwidth selection see \cite{bib5}.

This one-dimensional toy problem is nevertheless of interest in applications for constructing (PoS). For instance we can model an avenue as a one dimensional domain where predictions of the most dangerous places in a selected  zone can be  accomplished.

In the following sections we present the  non-conservation of the mean for the Neumann boundary conditions for  Problem (\ref{eq1}). We also propose new boundary conditions. For the derivations we assume that the functions are sufficiently smooth in order for the theorems of vector analysis to hold. Moreover the following derivations can be done for a more general diffusion equation with a variable diffusion coefficient $k(x)$.

\subsection{Neumann boundary conditions}

If we consider the Neumann or natural boundary conditions on the Problem (\ref{eq1}), we have

\begin{equation}\label{eq3}
 \frac{\partial u(x,t)}{\partial x}\Big|_{a}=0 \quad, \quad \frac{\partial u(x,t)}{\partial x}\Big|_{b}=0 .
\end{equation}

As is widely known, the total mass is conserved over time, see  Section \ref{mass-conv}, however the mean of the initial condition  is, in general, not conserved. Indeed, we have

\begin{align*}
\frac{d}{dt} \left( \int_{a}^{b} x u(x,t) dx \right) & =  \int_{a}^{b} x  \frac{\partial^2 u(x,t)}{\partial x^2} dx\\
& = \left[ x \frac{\partial u(x,t)}{\partial x} \right]_{a}^{b}- \left[u(x,t)\right]_{a}^{b}\\
& = u(a,t) - u(b,t).
\end{align*}

Where we used (\ref{eq1}), (\ref{eq3}) and integration by parts. Hence the mean is generally not conserved, it depends on the values of $u(x,t)$ at the boundary in a time $t$.

\subsection{Boundary conditions that conserve the mean}\label{mass-conv}

We propose the following boundary conditions for (\ref{eq1}),

\begin{equation}\label{eq4}
 \frac{\partial u(x,t)}{\partial x}\Big|_{a}=\frac{\partial u(x,t)}{\partial x}\Big|_{b} \quad, \quad \frac{u(b)-u(a)}{b-a} = \frac{\partial u(x,t)}{\partial x}\Big|_{b} .
\end{equation}

Note that this boundary conditions are non-local, we need to evaluate in  both boundary points at the same time. Now we show that both the mean and the mass are conserved over time using this boundary conditions. Consider first the conservation of the total mass. We have,

\begin{align*}
\frac{d}{dt} \left( \int_{a}^{b} u(x,t) dx \right)  =  \int_{a}^{b}   \frac{\partial^2 u(x,t)}{\partial x^2} dx
 = \left[  \frac{\partial u(x,t)}{\partial x} \right]_{a}^{b} = \frac{\partial u(x,t)}{\partial x}\Big|_{a}-\frac{\partial u(x,t)}{\partial x}\Big|_{b}=0.
\end{align*}

Where we used (\ref{eq1}), (\ref{eq4}) and integration by parts. This shows that the total mass is conserved. Consider now the conservation of the mean. We have,

\begin{align*}
\frac{d}{dt} \left( \int_{a}^{b} x u(x,t) dx \right) & =  \int_{a}^{b} x  \frac{\partial^2 u(x,t)}{\partial x^2} dx\\
& = \left[ x \frac{\partial u(x,t)}{\partial x} \right]_{a}^{b}- \left[u(x,t)\right]_{a}^{b}\\
& = (b-a)\frac{\partial u(x,t)}{\partial x}\Big|_{b} -u(b,t) + u(a,t)\\
& = 0.
\end{align*}

Again (\ref{eq1}), (\ref{eq4}) and integration by parts were used to obtain the desired result.

This shows that the boundary conditions (\ref{eq4}) for  problem (\ref{eq1}) conserve both mean and mass. Now we proceed to make some numerical simulations using FEM to show the consequences of the application of this boundary conditions in the process of estimation a probability density for a data sample (\ref{eq2}). 

\section{Numerical study of mean conservation }

Now the problem (\ref{eq1}),(\ref{eq4}) is written in a weak formulation \cite{bib7}  in order to apply the finite element method to the  problem. Now for all $v(x) \in C^{\infty}(a,b)$ we have,

\begin{equation} \label{eq5}
\int_{a}^{b} \frac{\partial u(x,t)}{\partial t} v(x)dx + \int_{a}^{b} \frac{\partial u(x,t)}{\partial x} \frac{d v(x)}{d x}dx =(v(b)-v(a))\frac{\partial u(x,t)}{\partial x}\Big|_{b}.
\end{equation}

We solve this weak formulation using FEM with low order elements in the\linebreak interval $[a,b]=[0,10]$, where the number of elements is $M$. Then  Problem (\ref{eq5}),(\ref{eq2}),(\ref{eq4}) yields the problem in the discretised space $V^{h}$. Find $u(x,t) \in V^{h}$, such that\linebreak for all $v(x) \in V^{h}$:

\begin{numcases}{}
   \int_{a}^{b} \frac{\partial u(x,t)}{\partial t} v(x)dx + \int_{a}^{b} \frac{\partial u(x,t)}{\partial x} \frac{d v(x)}{d x}dx =(v(b)-v(a))\frac{\partial u(x,t)}{\partial x}\Big|_{b}, \label{eq6}
   \\
    u(x,0)=\frac{M}{(b-a)N}\sum_{i=1}^{N}\delta(x-b_i), \ \ x,b_i  \in [a,b], \label{eq7}
      \\
    \frac{\partial u(x,t)}{\partial x}\Big|_{a}=\frac{\partial u(x,t)}{\partial x}\Big|_{b} \quad, \quad \frac{u(b)-u(a)}{b-a} = \frac{\partial u(x,t)}{\partial x}\Big|_{b}.  \label{eq8}
\end{numcases}

Where we represent delta measures by the closest base element of the finite element approximation. Note that (\ref{eq7}) contains a normalization factor, since now the elements integral are not one (since they are not delta measures).

Now we use the Galerkin method of mean weighted residuals for the spatial part of the problem choosing low order elements $\phi_{i}$. This formulation can be found in \cite{bib7}. For our numerical studies we solve the temporal part of the problem (element coefficients) using the implicit-Euler Galerkin Discretization \cite{bib8}, thus the problem is reduced to solve a linear system iteratively for every timestep $\Delta t$. 

\begin{figure}[h]
\begin{subfigure}{.47\textwidth}
  \centering
  \includegraphics[width=1\linewidth]{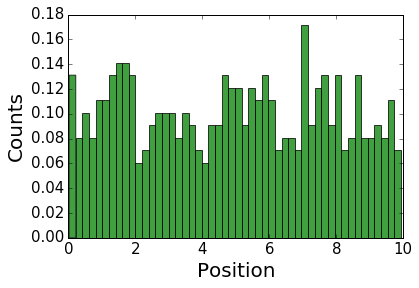}
  \caption{Histogram of the initial discrete data sample using 50 bins.}\label{histogram}
\end{subfigure}\hfill
\begin{subfigure}{.47\textwidth}
  \centering
  \includegraphics[width=1\linewidth]{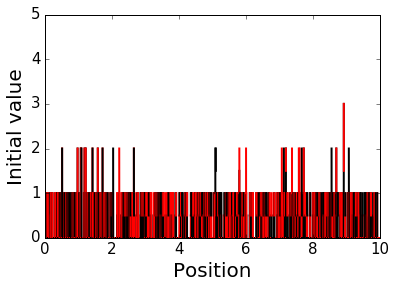}
  \caption{Initial discrete sample date seen as a plot of the FEM initial condition. }\label{initial-condition}
\end{subfigure}
\caption{Initial discrete data sample to estimate a continuous probability density, its chosen to be uniformly distributed in the interval $[0,10]$.}
\label{fig:fig1}
\end{figure}

In order to implement the previous formulation numerically, we use \texttt{python} to do all the calculations for the simulation. The code is available publicly in \cite{bib9}. There we start by generating a list of $\{b_n\}_{n=1}^{N=500}$   uniformly distributed points in the\linebreak interval  [0,10]. This points are located in the closest interval of the spatial FEM partition $\{ (0+ (n-1)/500,n/500) \}_{n=1}^{5000}$. The histogram of this points, Figure \ref{histogram}  can be seen for instance as the number of times a certain criminal act was informed in a zone from the avenue. See Figure \ref{fig:fig1}.

If we represent this data as an initial condition \eqref{eq7} we obtain the Figure \ref{initial-condition}. Where we plotted alternatively each consecutive FEM basis function red and black.

\begin{figure}[h!]
\begin{subfigure}{.47\textwidth}
  \centering
  \includegraphics[width=1\linewidth]{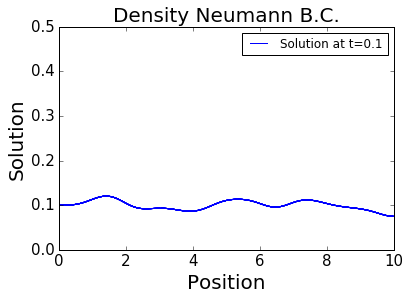}
  \caption{Numerical solution of the KDE problem  using Neumann boundary conditions  evolved a time $t=0.1$        .}\label{neumman-sol}
\end{subfigure}\hfill
\begin{subfigure}{.47\textwidth}
  \centering
  \includegraphics[width=1\linewidth]{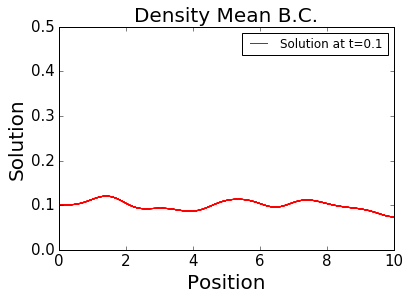}
  \caption{Numerical solution of the KDE problem  using mean conserving boundary conditions  evolved a time $t=0.1$ . }\label{mean-sol}
\end{subfigure}
\caption{Plots of the numerical solution of the problem \eqref{eq6},\eqref{eq7},\eqref{eq8} using both boundary conditions  evolved a time $t=0.1$ .}
\label{fig:fig2}
\end{figure}

Now  we solve numerically the problem using the implicit-Euler Galerkin discretization for the problem \eqref{eq6},\eqref{eq7},\eqref{eq8} and we evolve the solution until time $t=0.1$ using either Neumann boundary conditions, see Figure \ref{neumman-sol} and  mean conserving boundary conditions, see Figure \ref{mean-sol}. The solution for the mean conserving boundary condition is positive for this numerical experiment, see Figure \ref{mean-sol}, this fact  is currently being explored for future analytical studies.

\begin{figure}[h!]
\begin{subfigure}{.47\textwidth}
  \centering
  \includegraphics[width=1\linewidth]{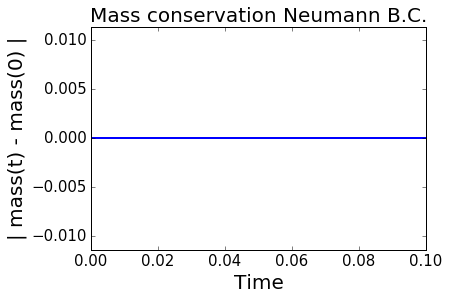}
  \caption{Change of mass $\Delta m$ for the numerical solution with Neumann boundary conditions. }
  \label{Mass-conv-Neu}
\end{subfigure}\hfill
\begin{subfigure}{.47\textwidth}
  \centering
  \includegraphics[width=1\linewidth]{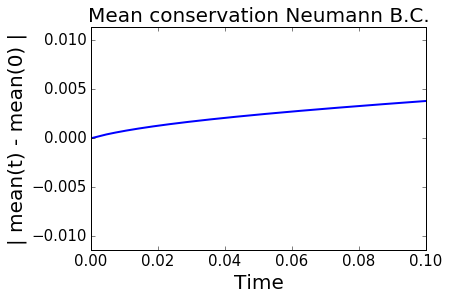}
  \caption{Change of mean $\Delta \mu$ for the numerical solution with Neumann boundary conditions. }
  \label{Mean-conv-Neu}
\end{subfigure}
\caption{Plots of the evolution of $\Delta m$ and $\Delta \mu$ for the density estimation with Neumann boundary conditions for $t \in [0,01]$.}
\label{fig:fig3}
\end{figure}

\begin{figure}[h]
\begin{subfigure}{.47\textwidth}
  \centering
  \includegraphics[width=1\linewidth]{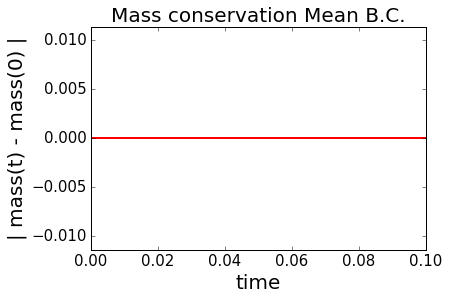}
  \caption{Change of mass $\Delta m$ for the numerical solution with mean conserving boundary conditions. }
  \label{Mass-conv-Mean}
\end{subfigure}\hfill
\begin{subfigure}{.47\textwidth}
  \centering
  \includegraphics[width=1\linewidth]{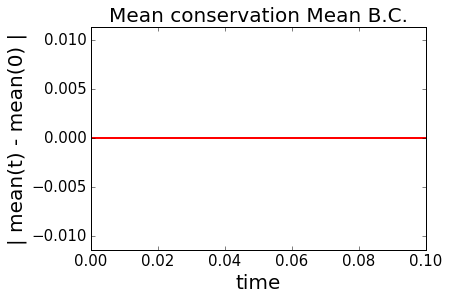}
  \caption{Change of mean $\Delta \mu$ for the numerical solution with mean conserving boundary conditions. }
  \label{Mean-conv-Mean}
\end{subfigure}
\caption{Plots of the evolution of $\Delta m$ and $\Delta \mu$ for the density estimation with mean conserving boundary conditions for $t \in [0,01]$.}
\label{fig:fig4}
\end{figure}

As the Figure \ref{fig:fig2} shows, the solutions are similar and therefore we can see that for this example the new boundary condition does not generate a noticeable change on the generation of the continuous density distribution. Nevertheless we present the plots of  change of mass $\Delta m(t) = m(0) - m(t)$, Figures \ref{Mass-conv-Neu}, \ref{Mass-conv-Mean} and change of mean $\Delta \mu(t) = \mu(0) - \mu(t)$, Figures \ref{Mean-conv-Neu}, \ref{Mean-conv-Mean} for both Neumann and mean conserving boundary conditions.

Figures \ref{Mean-conv-Neu} and \ref{Mean-conv-Mean} present the real difference in the evolution of the density. We effectively see that the mean conserving boundary conditions conserve the mean in the density estimation process. On the other hand if we where to have an initial condition that is biased to one of the boundaries, the differences of the estimated densities by both boundary conditions would differ significantly. However there is no evidence to think that this phenomena occurs in real avenues.

For the numerical experiment presented here we can see that the mean for the Neumann boundary conditions has changed about 0.4\% in $t = 0.1$. This change is small, in fact, for an avenue of 10 km, the change in mean would be about 40 m.  We conclude that for this numerical experiment  for the process of density estimation (when the data has not change to much due to the smoothing process) the Neumann boundary condition provide a very fast (since they are easy to implement) and  accurate way to estimate a continuous probability density.  Nevertheless the mean of the sample is not preserved exactly, on the other hand, the mean conserving boundary condition, apart from being also easily implementable, is accurate and do preserve the mean of the sample.

\section{Two-dimensional densities}

We now present the problem for the diffusion equation \cite{bib10} in two dimensions,

\begin{equation}\label{eq9}
\frac{\partial u({\bf x},t)}{\partial t}- \nabla^{2} u({\bf x},t) =0, \qquad {\bf x}=(x_1,x_2) \in \Omega \subset \mathbb{R}^2 , t>0.
\end{equation}

Again we want the conservation of mass and mean in the time evolution of the density. Consider first the conservation of the total mass. We have,

\begin{align*}
\frac{d}{dt} \left( \int_{\Omega} u({\bf x},t) d{\bf x} \right)  =  \int_{\Omega}   \nabla^{2} u({\bf x},t)  d{\bf x}
 = \int_{\partial \Omega}    \frac{\partial u({\bf x},t)}{\partial \boldsymbol{\nu}}  d{ \sigma},
\end{align*}

where $\nabla u \cdot \boldsymbol{\nu}= \frac{\partial u({\bf x},t)}{\partial \boldsymbol{\nu}} $, and $\boldsymbol{\nu}$ denotes  the outward normal unit vector to $\partial \Omega$.   To deduce this relation we used (\ref{eq9}), and the first Green identity \cite{bib10}. Consider now the conservation of the mean. We have,

\begin{align*}
\frac{d}{dt} \left( \int_{\Omega} x_i u({\bf x},t) d{\bf x} \right) & =   \int_{\Omega}  x_i \nabla^{2} u({\bf x},t)  d{\bf x}\\
& =\int_{\partial \Omega}     x_i \frac{\partial u({\bf x},t)}{\partial \boldsymbol{\nu}}  d{ \sigma} - \int_{\Omega}  \nabla_i u({\bf x},t) d{\bf x},
\end{align*}

where $\nabla_i u({\bf x},t) =  {\bf e}_i  \cdot \nabla u({\bf x},t)$, assuming Cartesian unit vectors. Again (\ref{eq9}) and the first Green's identity   were used to obtain the desired result.

Then the conditions that we have to impose on $u({\bf x},t)$ in order to conserve mean and mass are:

\begin{equation}
\int_{\partial \Omega}     x_i \frac{\partial u({\bf x},t)}{\partial \boldsymbol{\nu}}  d{ \sigma} = \int_{\Omega}  \nabla_i u({\bf x},t) d{\bf x}  , \quad i = 1,2, \qquad \text{and} \qquad \int_{\partial \Omega}    \frac{\partial u({\bf x},t)}{\partial \boldsymbol{\nu}}  d{ \sigma} =0.
\end{equation}

The advantage of two dimensional domains is that we are not restricted to impose only two conditions for the boundary(mean and mass conservation). For these domains we can in principle conserve additional higher moments of the density distribution that are meaningful for the particular problem. Applications on two dimensional domains are of special interest for the author since a two dimensional map of the city can generate really robust results in the field of  Perception of security(PoS). 

\section{Conclusions}

The proposed mean conserving boundary conditions were shown to effectively maintain the mean of the initial data sample over the continuous density estimation process. This was also confirmed by the numerical simulation of the estimation process where we used a list of uniformly distributed points in the interval [0,10] as an initial condition.

The numerical experiments presented here show that even though Neumann boundary conditions do not conserve the mean over time, they are accurate enough to maintain the mean in a very restricted interval before the over-smoothing of the density estimation process. 

We showed the application and some of the consequences of both the idea of (KDE) and the new boundary conditions to avenues in a city. The consequences of implementing the diffusion equation with the proposed boundary conditions in companion of more special initial conditions and in 2D domains remains to be analyzed. 

\section*{Acknowledgments}

I would like to express my gratitude to  Juan Galvis, whose guidance was essential for the completion of this manuscript. I also want to thank Francisco A. G\'omez and Zdravko Botev, whose comments were really appreciated for the analysis of the results.

\end{document}